\documentclass[a4paper]{amsart}

\usepackage{amsmath}
\usepackage{amscd}
\usepackage{fullpage}

\usepackage[utf8]{inputenc}
\usepackage[]{geometry}

\usepackage[utf8]{inputenc}
\usepackage[T1]{fontenc}
\usepackage{lmodern}

\usepackage{listings}
\makeatletter
\def\@fnsymbol#1{\ensuremath{\ifcase#1\or \star\or\dagger\or \ddagger\or
   \mathsection\or \mathparagraph\or \|\or **\or \dagger\dagger
   \or \ddagger\ddagger \else\@ctrerr\fi}}
\makeatother

\usepackage{amsmath}
\usepackage{mathtools}
\usepackage{amssymb}
\usepackage{amsthm}   

\usepackage{graphicx}

\usepackage{tikz}
\usepackage{ifthen}
\usetikzlibrary{math}
\usetikzlibrary{fit}

\newcommand{\BeHe}[6]{
\begin{tikzpicture}[scale=1/#1*15,cap=round,>=latex]
\pgfmathparse{0};
\edef\shiftvar{\pgfmathresult};

\foreach \n in {1,...,#1}{

    \pgfmathparse{gcd(\n,#1)};
    
    \ifthenelse{1=\pgfmathresult}{
    \draw[thick,xshift=3*\shiftvar cm] (0cm,0cm) circle(1cm);
    \filldraw[red,xshift=3*\shiftvar cm] (\n*#2*360:1cm) circle(3pt);
    \filldraw[red,xshift=3*\shiftvar cm] (\n*#3*360:1cm) circle(3pt);
    \filldraw[red,xshift=3*\shiftvar cm] (\n*#4*360:1cm) circle(3pt);
    \filldraw[blue,xshift=3*\shiftvar cm] (\n*#5*360:1cm) circle(3pt);
    \filldraw[blue,xshift=3*\shiftvar cm] (\n*#6*360:1cm) circle(3pt);
    \filldraw[blue,xshift=3*\shiftvar cm] (\n*1*360:1cm) circle(3pt);
    \pgfmathparse{\shiftvar+1} \xdef\shiftvar{\pgfmathresult}
    }{};
}
\end{tikzpicture}
}

\usepackage{array}

\usepackage[plainpages=false,pdfpagelabels,colorlinks=true,citecolor=blue,hypertexnames=false]{hyperref}

\newtheorem{coro}{Corollary}
\newtheorem{thm}{Theorem}

\newtheorem{lem}{Lemma}

\theoremstyle{definition}

\newcommand{\ps}[1]{[\![#1]\!]}
\newcommand{\NN}[0]{\mathbb{N}}
\newcommand{\ZZ}[0]{\mathbb{Z}}
\newcommand{\QQ}[0]{\mathbb{Q}}

\newcommand{\x}{\mathbf{x}}
\newcommand{\K}{\mathbb{K}}

\newcommand{\diag}[0]{\mathrm{Diag}}

\usepackage[textwidth=1.7cm]{todonotes}

\protected\def\hrefalin{\href{https://specfun.inria.fr/bostan/}}

\protected\def\hrefsergey{\href{https://homepage.univie.ac.at/sergey.yurkevich/}}

\begin{document}

 \title{On a Class of Hypergeometric Diagonals}
\author[Alin Bostan]{\hrefalin{Alin Bostan}}
\address{{\href{{https://specfun.inria.fr/bostan/}}{Alin Bostan}}, Inria, Universit\'e Paris-Saclay, 
1 rue Honor\'e d'Estienne d'Orves, 91120,
Palaiseau, France}
\email{\href{mailto:alin.bostan@inria.fr}{\tt alin.bostan@inria.fr}}
\thanks{Alin Bostan was supported in part by \textcolor{magenta}{\href{https://specfun.inria.fr/chyzak/DeRerumNatura/}{DeRerumNatura}} ANR-19-CE40-0018.}

\author[Sergey Yurkevich]{\hrefsergey{Sergey Yurkevich}}
\address{{\href{{https://homepage.univie.ac.at/sergey.yurkevich/}}{Sergey Yurkevich}}, University of Vienna, Austria}
\email{\href{mailto: sergey.yurkevich@univie.ac.at}{\tt sergey.yurkevich@univie.ac.at}}
\thanks{Sergey Yurkevich was supported by the \href{https://www.fwf.ac.at/}{Austrian Science Fund} (FWF): P-31338.}
                                     



                       
\dedicatory{}

\begin{abstract}
We prove that the diagonal of any finite product of algebraic functions of 
the form
\begin{align*} 
    {(1-x_1- \dots -x_n)^R},
    \qquad R\in\mathbb{Q},
\end{align*}
is a generalized hypergeometric function, and we provide an explicit description of its parameters. The particular case $(1-x-y)^R/(1-x-y-z)$
corresponds to the main identity of Abdelaziz, Koutschan and
Maillard in~\cite[\S3.2]{AbKoMa20}. Our result is useful in both directions: on
the one hand it shows that Christol's conjecture holds true for a large class
of hypergeometric functions, on the other hand it allows for a very explicit
and general viewpoint on the diagonals of algebraic functions of the type
above. Finally, in contrast to~\cite{AbKoMa20}, our proof is completely
elementary and does not require any algorithmic help.  
\end{abstract}         

\maketitle


\section{Introduction}

Let $\K$ be a field of characteristic zero and 
let $g \in \K\ps{\x}$ be a power series in $\x = (x_1,\dots,x_n)$  
\begin{align*}
    g(\x) &= \sum_{(i_1,\dots,i_n)\in\NN^n} g_{i_1,\dots,i_n}
    x^{i_1}_1 \cdots x^{i_n}_n \in \K[\![\x]\!].
\end{align*}   
The \emph{diagonal} $\diag(g)$ of $g(\x)$ is the univariate power series given by
\begin{align*}
    \diag(g)
    \coloneqq \sum_{j \geq 0} g_{j,\dots,j} t^j \in \K[\![t]\!].
\end{align*}
A power series $h(\x)$ in $\K[\![\x]\!]$ is called \emph{algebraic} if there exists a non-zero
polynomial $P(\x,T) \in \K[\x,T]$ such that $P(\x,h(\x)) =0$; otherwise, it is
called \emph{transcendental}.           

If $g(\x)$ is algebraic, then its diagonal $\diag(g)$ is usually
transcendental; however, by a classical result by Lipshitz~\cite{Lipshitz88},
$\diag(g)$ is D-finite, i.e., it satisfies a non-trivial linear differential
equation with polynomial coefficients in $\K[t]$. Equivalently, the coefficients
sequence $(g_{j,\dots,j})_{j \geq 0}$ of $\diag(g)$ is P-recursive, i.e., it
satisfies a linear recurrence with polynomial coefficients (w.r.t.~$j$).

When a P-recursive sequence satisfies a recurrence of order 1, we say that it is \emph{hypergeometric}.
An important class of power series, whose coefficients sequence is hypergeometric by design, is that of {generalized hypergeometric functions}.
Let $p,q \in \NN$ and $a_1,\dots,a_p$ and $b_1,\dots,b_q$ be rational numbers such that $b_i+j \neq 0$ for any $i,j \in \NN$. The \emph{generalized hypergeometric function} ${}_{p}F_{q}$ with parameters $a_1,\dots,a_p$ and $b_1,\dots,b_q$ is the univariate power series in $\K[\![t]\!]$ defined by
\[
\,{}_{p}F_{q}([a_{1},\ldots ,a_{p}]\,;[b_{1},\ldots ,b_{q}]\,;t) \coloneqq \sum _{j\geq 0}{\frac {(a_{1})_{j}\cdots (a_{p})_{j}}{(b_{1})_{j}\cdots (b_{q})_{j}}}\,{\frac {t^{j}}{j!}},
\]
where $(x)_j \coloneqq x(x+1)\cdots(x+j-1)$ is the rising factorial.

\medskip 
We are interested in this article in the following (dual) questions: 
\begin{itemize}
	\item[(i)] What are the algebraic power series $g(\x)$ whose diagonal $\diag(g)$ is a generalized hypergeometric function ${}_{p}F_{q}$?\footnote{Note that a necessary condition is that $q=p-1$, since the radius of convergence must be finite and non-zero.} 
	\item[(ii)] What are the hypergeometric sequences $(a_j)_{j \geq 0}$ whose generating functions $\sum_{j \geq 0} a_j t^j$ can be written as diagonals of algebraic power series? 
\end{itemize}	

Already for $n\in \{ 1, 2 \}$ these questions\footnote{From an algorithmic viewpoint, questions (i) and (ii) are very different in nature: while (i) is decidable (given an algebraic power series, one can decide if its diagonal is hypergeometric, for instance by combining the algorithms in~\cite{BoLaSa13} and~\cite{Petkovsek92}), the status of question (ii) is not known (does there exist an algorithm which takes as input a hypergeometric sequence and outputs an algebraic series whose diagonal is the generating function of the input sequence?).} are non-trivial. The classes of
diagonals of bivariate rational power series and of algebraic power series
coincide~\cite{Polya22,Furstenberg67}. Hence, questions (i) and (ii) contain
as a sub-question the characterization of algebraic hypergeometric functions.
This problem was only recently solved in a famous paper by Beukers and
Heckman~\cite{BeHe89}.

\smallskip
Another motivation for studying questions (i) and (ii) comes from the well-known conjecture below, formulated in~\cite{Christol87,Christol90} by Christol. Recall that $f \in \QQ\ps{t}$ is called \emph{globally bounded} if it has finite non-zero radius of convergence and $\beta \cdot f (\alpha \cdot t) \in \ZZ\ps{t}$ for some non-zero $\alpha, \beta\in\ZZ$. 
\begin{quote}
{\bf Christol's conjecture.} If  $f\in \mathbb{Q}\ps{t}$ is
D-finite and globally bounded, then $f=\diag(g)$ for some $n\in\NN$ and some
algebraic power series $g\in \mathbb{Q}[\![x_1, \ldots, x_n]\!]$. 
\end{quote}

Christol's conjecture is still largely open, even in the particular case when
$f$ is a generalized hypergeometric function. In this case, it has
been proved~\cite{Christol87,Christol90} in two extreme subcases: when all the bottom parameters $b_i$ are
integers (case of ``minimal monodromy weight'', in the terminology of~\cite{Christol15}) and when they are all non-integers
(case of ``maximal monodromy weight''). 
In the first extremal case, the proof is based on the  observation that
\begin{equation}
	_{p}F_{q}([a_{1},\ldots ,a_{p}]\,;[{1},\ldots ,{1}]\,;t)
	=
	(1-t)^{-a_1} \star \cdots \star (1-t)^{-a_p},
\end{equation}                                     
where $\star$ denotes the Hadamard (term-wise) product, and on the fact that
diagonals are closed under Hadamard product~\cite[Prop.~2.6]{Christol88}.     
In the second extremal case,
it is based on the equivalence between being globally bounded and algebraic; this equivalence, proved by Christol~\cite{Christol87,Christol90}, is itself based on~\cite{BeHe89}.     
                  
The other cases (of ``intermediate monodromy weight'') are widely open.
A first explicit example of this kind, itself still open as of today, was given by Christol himself as soon as 1987~\cite[\S VII]{Christol87}:
\begin{quote}  
Is $f(t) = {}_{3}F_{2} \left( \left[ \frac19, \frac49, \frac59 \right]\,; \left[1, \frac13 \right]\,;t \right)$
the diagonal of an algebraic power series?
\end{quote}  

Two decades later, Bostan et al.~\cite{BBCHM12,BBCHM13} produced a large list of about 100 similar $_{3}F_{2}$ (globally bounded) functions, which are potential counter-examples to Christol's conjecture (in the sense that, like ${}_{3}F_{2}([1/9,4/9,5/9]\,;[1,1/3]\,;t)$, they are not easily reducible to the two known extreme cases, via closure properties of diagonals, e.g., with respect to Hadamard products).
In 2020, Abdelaziz, Koutschan and Maillard~\cite[\S3]{AbKoMa20} managed to show that two members of that list, namely 
$_{3}F_{2}([1/9,4/9,7/9]\,;[1,1/3]\,;t)$
and
$_{3}F_{2}([2/9,5/9,8/9]\,;[1,2/3]\,;t)$
are indeed diagonals. Precisely,
\begin{equation} \label{eq:10}  
 {} _{3}F_{2} \left(\left[\frac29, \frac59, \frac89 \right]\,; \left[1, \frac23 \right]\,; 27 \, t  \right)
=
\mathrm{Diag} \left( \frac{(1-x-y)^{1/3}}{1-x-y-z} \right) 
\end{equation}                                 
and 
\begin{equation} \label{eq:11}  
{} _{3}F_{2} \left( \left[ \frac19, \frac49, \frac79 \right]\,; \left[1, \frac13 \right]\,; 27 \, t \right)
= 
\mathrm{Diag} \left( \frac{(1-x-y)^{2/3}}{1-x-y-z} \right). 
\end{equation}

Sec.~3 of \cite{AbKoMa20} also contains the following extension of identities~\eqref{eq:10} and~\eqref{eq:11}, to any $R\in\mathbb{Q}$:
\begin{equation} \label{eq:23-26}  
{} _{3}F_{2} \left( \left[ 
\frac{1-R}{3},\frac{2-R}{3},\frac{3-R}{3}
\right]\,;
[1, 1-R]\,; 27 \, t \right)
= 
\mathrm{Diag} \left( \frac{(1-x-y)^R}{1-x-y-z} \right). 
\end{equation}

A common feature of identities~\eqref{eq:10} and~\eqref{eq:11} (and 
their extension~\eqref{eq:23-26}) is that the top parameters are in
arithmetic progression, as opposed to Christol's initial example. However,
they are the first known examples of generalized hypergeometric functions with
intermediate monodromy weight, not trivially reducible to the two known extreme cases,
and which are provably diagonals.

\smallskip Our first result extends identity~\eqref{eq:23-26} to a much larger
class of (transcendental) generalized hypergeometric functions.


\begin{thm} \label{thm:1}
    Let $R,S \in \QQ$ and $n,N \in \NN$ such that $S \neq 0$ and $0 \leq n \leq  N$. Set $s \coloneqq N-n$ and $Q \coloneqq S-R$. Then the generalized hypergeometric function
    \begin{align*}
         & {}_{N+s} F_{N+s-1} \left(\left[\frac{Q}{N},\frac{Q+1}{N},\dots,\frac{Q+N-1}{N}, \frac{S}{s}, \dots, \frac{S+s-1}{s}\right]\,;
         \left[\frac{Q}{s}, \dots, \frac{Q+s-1}{s}, 1, \dots, 1\right]\,; N^Nt\right)
    \end{align*}
is equal to the diagonal    
\begin{align*}
         & \mathrm{Diag} \left( \frac{(1-x_1- \dots -x_n)^R}{(1-x_1 -\dots-x_N)^S} \right).
    \end{align*}
\end{thm}     
Note that identity~\eqref{eq:23-26} corresponds to the particular case
$(n,N,S) = (2,3,1)$ of Theorem~\ref{thm:1}. The proof of~\eqref{eq:23-26} given
in~\cite[\S3.2]{AbKoMa20} relies on an algorithmic technique called
\emph{creative telescoping}~\cite{Koutschan13}, which works in principle\footnote{Creative telescoping algorithms, such as the one in~\cite{BoLaSa13}, compute a linear differential equation for ${\rm Diag}(g(\x))$. This equation is converted to a linear recurrence, whose hypergeometric solutions can be computed using Petkov\v{s}ek's algorithm~\cite{Petkovsek92}.
Note that the
complexity (in time and space) of these algorithms increases with $n,N,R$ and
$S$.} on any diagonal of an algebraic function, as long as the number  $\max(n,N)$ of
indeterminates is \emph{fixed}. 
Our identity in
Theorem~\ref{thm:1} contains a number of indeterminates which is itself {variable}, hence
it cannot be proved by creative telescoping in this generality. In \S\ref{proof:thm3} we offer instead a direct and elementary proof of a natural generalization. We note that Theorem~\ref{thm:1} can also be proven by directly multiplying out the argument of the diagonal using the multinomial theorem, collecting needed terms and simplifying using the Chu-Vandermonde identity.

Note that in the theorem above, and similarly in later statements, the restriction on $R$ and~$S$ to be rational numbers is actually superfluous. Indeed, from the proofs it is obvious that the identities hold for arbitrary (formal) parameters $R,S$;  however, we include this condition because we wish $(1-x_1- \dots -x_n)^R(1-x_1 -\dots-x_N)^{-S}$ to be an algebraic function.

In \S\ref{sec:gen-Thm1} we will further generalize Theorem~\ref{thm:1}
in two distinct directions. The first extension (Theorem \ref{thm:general1})
shows that the diagonal of the product of an arbitrary number of arbitrary powers
of linear forms of the type $1 - x_1 - \cdots - x_m$ is again a generalized
hypergeometric function. The second extension (Theorem \ref{thm:general2}) shows that under a condition on the exponents the same stays true if the product is multiplied with another factor of the form $(1 - x_1 - \cdots - x_{m-2} - 2 \, x_{m-1})^b$.
 For instance, when restricted to $m=3$ variables, these results specialize as follows:

\begin{thm} \label{cor:23}
For any $R,S,T \in \mathbb{Q}$, we have:
\begin{align} \label{id:gen1}
    &\diag \left( (1-x)^R(1-x-y)^S(1-x-y-z)^T \right) =  \\
    &{}_6F_{5}\bigg( \left[ \frac{-(R+S+T)}{3},\frac{1-(R+S+T)}{3},\frac{2-(R+S+T)}{3}, \frac{-(S+T)}{2} , \frac{1-(S+T)}{2} , -T \right] \,; \nonumber \\ & \hspace{3cm} \left[  \frac{-(R+S+T)}{2},\frac{1-(R+S+T)}{2}, -(S+T), 1,1\right] \,; 27 t\bigg)  \nonumber  
\end{align}        
and
\begin{align}  \label{id:gen2}  
    \diag \left( (1-x)^R(1-x-2\, y)^S(1-x-y-z)^{-1} \right) & = \\
    {}_4F_{3}\bigg( \left[ \frac{1-(R+S)}{3},\frac{2-(R+S)}{3},\frac{3-(R+S)}{3}, \frac{1-S}{2} \right] & \,; \left[  \frac{1-(R+S)}{2},\frac{2-(R+S)}{2}, 1\right] \,; 27 t\bigg). \nonumber    
\end{align}      
\end{thm}

Note that~\eqref{id:gen2} generalizes and explains the  following two identities observed in~\cite[Eq.~(30)--(31)]{AbKoMa20}
\begin{equation} \label{eq:30}  
 {} _{3}F_{2} \left( \left[\frac19, \frac49, \frac79 \right]\,; \left[1, \frac23\right]\,; 27 \, t \right)
=
\mathrm{Diag} \left( \frac{(1-x- 2\, y)^{2/3}}{1-x-y-z} \right) 
\end{equation}                                 
and 
\begin{equation} \label{eq:31}   
	 {} _{3}F_{2} \left( \left[\frac29, \frac59, \frac89 \right]\,; \left[1, \frac56\right]\,; 27 \, t \right)
= 
\mathrm{Diag} \left( \frac{(1-x- 2\, y)^{1/3}}{1-x-y-z} \right). 
\end{equation}      
Once again, our proofs  of the (generalizations of) identities~\eqref{id:gen1} 
and~\eqref{id:gen2} are elementary, and do not rely on algorithmic tools.
   
One may wonder if other generalizations are possible, for instance whether the coefficient 2 can be replaced by a different one in~\eqref{id:gen2}.
The following example shows that this is  not the case.
Let
\[        
U(t) = \mathrm{Diag} \left( {\frac {\sqrt [3]{1-ax}}{1-x-y}} \right)  .
\] 
Then, the coefficients sequence $(u_j)_{j \geq 0}$ of $U(t)$
satisfies the second-order recurrence relation
\begin{align*}
 2  {a}^{2}  \left( 6 n+5 \right)  \left( 3 n+1 \right) u_n - 3 & \left( n+1 \right)  \left( 3 \left( {a}^{2}+4a-4 \right) n+2 {a}^{2}+18 a-18 \right) u_{n+1}     \\
 =  9 & \left( 1-a \right) \left( n+2 \right)  \left( n+1 \right)  u_{n+2}.
\end{align*}   
When $a\in \{ 0, 1, 2\}$, the sequence $(u_j)_{j \geq 0}$ also satisfies
a shorter recurrence, of order 1, as shown by our main results.
In these cases, $U(t)$ is a hypergeometric function. 
When $a\notin \{ 0, 1, 2\}$, the second-order recurrence is the minimal-order satisfied by $(u_j)_{j \geq 0}$, hence   $U(t)$ is not a hypergeometric function.
This can be proved either using the explicit identity
\[
U(t) = 
\sqrt [3]{{\frac {a/2}{1-4\,t}}+{\frac {1-a/2}{ \left( 1-4\,t \right) ^{
{\frac{3}{2}}}}}},
\]   
or by using the general approach in~\cite[\S5]{Boucher99}.

An apparent weakness of our results is that they only provide
examples with parameters in (unions of) arithmetic progressions.
This is true, as long as identities are used alone. But symmetries
may be broken by combining different identities and using for instance Hadamard products. As an illustration, by taking the Hadamard product in both sides of the following identities
\[
 {} _{3}F_{2} \left(\left[\frac{Q}{3}, \frac{Q+1}{3}, \frac{Q+2}{3} \right]\,; \left[1, Q \right]\,;  t  \right)
=
\mathrm{Diag} \left( \frac{(1-\frac{x_1}{3}-\frac{x_2}{3})^{1-Q}}{1-\frac{x_1}{3}-\frac{x_2}{3}-\frac{x_3}{3}} \right) 
\]
and
\[
 {} _{2}F_{1} \left(\left[\frac{Q}{6}, \frac{Q+3}{6} \right]\,; \left[ \frac{Q}{3} \right]\,; t  \right)
=
\mathrm{Diag} \left( \frac{(1-\frac{x_4}{2})^{1-Q/3}}{1-\frac{x_4}{2}-\frac{x_5}{2}} \right), \]
both particular cases of Theorem ~\ref{thm:1}, one deduces that the \emph{non-symmetric} hypergeometric function
\[
{
{}_4F_{3}
\left(
\left[{\frac {Q}{6}},{\frac {Q+3}{6}},{\frac {Q+1}{3}},{\frac {Q+2}{3}}\right]\,;
\left[1,1,Q\right]\,;t
\right)}
\]      
is equal to the diagonal 
\[    
\mathrm{Diag} \left(
{ \left( 1-{\frac {x_1}{3}}-{\frac {x_2}{3}} \right) ^{1-Q}
 \left( 1-{\frac {x_4}{2}} \right) ^{1-{\frac {Q}{3}}} \left( 1-{
\frac {x_1}{3}}-{\frac {x_2}{3}}-{\frac {x_3}{3}}
 \right) ^{-1} \left( 1-{\frac {x_4}{2}}-{\frac {x_5}{2}}
 \right) ^{-1}}
\right).
\]

Similarly, the \emph{non-symmetric} hypergeometric function\footnote{Amusingly, the above $_3F_{2}$ is not only asymmetric, but it also shares another similarity with  Christol's example: the sum of two of the three top parameters is equal to the third one. This pattern occurs in several other examples.}               
\[
{
{}_3F_{2}
\left(
\left[{\frac {Q}{6}},{\frac {Q+3}{6}},{\frac {2Q+3}{6}}\right]\,; 
\left[1,{\frac {2\,Q}{3}}\right]\,;t
\right)}
\]  
is equal to the diagonal 
\[
\mathrm{Diag}  \left( { \left( 1-{\frac {x_1}{2}} \right) ^{1-{
\frac {Q}{3}}} \left( 1-{\frac {x_3}{2}} \right) ^{1-{\frac {2\,Q
}{3}}} \left( 1-{\frac {x_1}{2}}-{\frac {x_2}{2}} \right) ^{
-1} \left( 1-{\frac {x_3}{2}}-{\frac {x_4}{2}} \right) ^{-1}
} \right)   .
\]
A natural challenge is to prove (or disprove) that Christol's 
${}_3F_{2}$ can be obtained in such a way.     

\smallskip

As a final remark,
one should not think that every generalized hypergeometric function which is 
a diagonal needs to have a representation like in our Theorems~\ref{thm:1} or~\ref{cor:23}.
For instance, the diagonal $\mathrm{Diag} \left((1-(1+w)(x+y+z))^{-1} \right) $          
is equal \cite{BBMW15} to the generalized hypergeometric function
\[
{
{}_4F_{3}
\left(
\left[{\frac {1}{3}}, \frac13, \frac23, {\frac{2}{3}}\right]\,; 
\left[1,1, {\frac {1}{2}}\right]\,;\frac{729}{4} \,t
\right)} 
=
1+18\,t+1350\,{t}^{2}
+ \cdots,
\]  
which is seemingly not of the form covered by any of our results.


\section{General case} \label{sec:gen-Thm1}
This section contains several parts: first we introduce in \S\ref{subsec:gen-Thm1} and \S\ref{subsec:gen-Thm1b} some notation and state the two general Theorems \ref{thm:general1} and \ref{thm:general2}. Then we explain them in~\S\ref{sec:examples} by means of four examples, showing that both Theorem~\ref{thm:1} and Theorem~\ref{cor:23} are special cases. Further, we continue in \S\ref{sec:lemmasproofs} with several lemmas and their proofs. Finally, the general theorems are proven in \S\ref{proof:thm3} and \S\ref{proof:thm4}.

\subsection{First Statement}\label{subsec:gen-Thm1}
Let $N \in \NN\setminus\{ 0 \}$ and $b_1,\dots,b_N \in \QQ$ with $b_N \neq 0$. We want to prove that the diagonal of
\begin{align} \label{eq:genR}
    R(x_1,\dots,x_N) \coloneqq (1+x_1)^{b_1} (1+x_1+x_2)^{b_2} \cdots (1+x_1+\cdots + x_N)^{b_N} 
\end{align}
can be expressed as a hypergeometric function. For each $k= 1,\dots,N$ we define the tuple 
\[
u^k \coloneqq \left( \frac{B(k)}{N-k+1}, \frac{B(k)+1}{N-k+1} , \dots , \frac{B(k)+N-k}{N-k+1}  \right),
\]
where $B(k) \coloneqq -(b_k+\cdots+b_N)$. For $k= 1, \dots ,N-1$ we set
\[
v^{k} \coloneqq \left( \frac{B(k)}{N-k}, \frac{B(k)+1}{N-k} , \dots , \frac{B(k)+N-k-1}{N-k}  \right).
\]
Moreover set $v^N \coloneqq (1,1,\dots,1)$ with exactly $N-1$ ones. It follows by construction that the lengths of the tuples 
\begin{align*}
    u & \coloneqq (u^1,\dots,u^N) \quad \text{ and } \quad v \coloneqq (v^1,\dots,v^N)
\end{align*}
are given by $M \coloneqq N+\cdots+2+1 = N(N+1)/2$ and $M-1$ respectively. We have the following generalization of Theorem~\ref{thm:1}:
\begin{thm} \label{thm:general1}
    It holds that
    \begin{align*}
        \diag(R(x_1,\dots,x_N)) = {}_{M}F_{M-1}([u]\,;[v]\,;(-N)^N t).
    \end{align*}
\end{thm}

\subsection{Second Statement}\label{subsec:gen-Thm1b}        
Let $N \in \NN\setminus\{ 0 \}$ and $b_1,\dots,b_N \in \QQ$. Assume that $b_{N} \neq 0$ and $b_{N-1}+b_{N} = -1$. We will prove that,
for any $b \in \QQ$,  we can express 
\begin{align*}
    (1+x_1+\cdots+x_{N-2}+2\, x_{N-1})^{b} \cdot R(x_1,\dots,x_N) 
\end{align*}
as a hypergeometric function as well. Again, let $B(k) \coloneqq -(b_k+\cdots+b_N)$. For each $k= 1,\dots,N-2$ we define the tuple 
\[
\tilde{u}^k \coloneqq \left( \frac{B(k) - b}{N-k+1}, \frac{B(k) - b+1}{N-k+1} , \dots , \frac{B(k) - b +N-k}{N-k+1}  \right)
\]
and set $\tilde{u}^{N-1} \coloneqq -(b_{N-1}+b_N+b)/2 = (1-b)/2$ and $\tilde{u}^N \coloneqq -b_N$. Moreover, for $k= 1, \dots ,N-2$ we set 
\[
\tilde{v}^{k} \coloneqq \left( \frac{B(k) - b}{N-k}, \frac{B(k) - b+1}{N-k} , \dots , \frac{B(k) - b+N-k-1}{N-k}  \right),
\]
and $\tilde{v}^{N-1} \coloneqq (1,1,\dots,1)$ with exactly $N-1$ ones. It follows by construction that the lengths of the tuples 
\begin{align*}
    \tilde{u} & \coloneqq (\tilde{u}^1,\dots,\tilde{u}^N) \quad \text{ and } \quad \tilde{v} \coloneqq (\tilde{v}^1,\dots,\tilde{v}^{N-1})
\end{align*}
are given by $M -1 =  N+\cdots+4+3 + 1+ 1 = N(N+1)/2-1$ and $M-2$ respectively.
\begin{thm}\label{thm:general2}
    It holds that
    \begin{align*}
        \diag((1+x_1+\cdots+x_{N-2}+2\, x_{N-1})^{b} \cdot R(x_1,\dots,x_N) ) = {}_{M-1}F_{M-2}([\tilde{u}]\,;[\tilde{v}]\,;(-N)^N t).
    \end{align*}
\end{thm}

\subsection{Examples} \label{sec:examples}
Let us list some examples of the general theorems and draw the connection to previous statements.
\begin{enumerate}
\item First we emphasize that Theorem \ref{thm:1} follows promptly from the more general Theorem \ref{thm:general1} by letting all $b_j = 0$ except $b_n = R$ and $b_N = -S$. Clearly, the change $\textbf{x} \mapsto -\textbf{x}$ in the algebraic function is reflected by the change $t \mapsto (-1)^N t$ in its diagonal.
\item Letting $N=3$ in Theorem \ref{thm:general1} we obtain immediately the first part of Theorem \ref{cor:23}. 
\item If moreover $T=-1$ in the case $N=3$, we achieve a cancellation of the last parameter and are left with
\begin{align*}
    &\diag \left( \frac{(1+x)^R(1+x+y)^S}{1+x+y+z} \right) = \\
    & \hspace{1cm} {}_5F_{4}\bigg( \left[ \frac{1-(R+S)}{3},\frac{2-(R+S)}{3},\frac{3-(R+S)}{3}, \frac{1-S}{2} , \frac{2-S}{2} \right] \,; \\ 
    & \hspace{3cm}  \left[  \frac{1-(R+S)}{2},\frac{2-(R+S)}{2}, 1-S, 1\right] \,; -27 t\bigg).    
\end{align*}
\item Comparing with the similar situation of Theorem \ref{thm:general2} in the case $N=3$ and  $b_{N-1} = -1 - b_N = 0$, we see that a family of ${}_4F_{3}$ functions remains and covers the second statement of Theorem \ref{cor:23}. 
\end{enumerate}

\subsection{Lemmas and Proofs}\label{sec:lemmasproofs}
In this section we will state and prove necessary lemmas for the proofs of Theorems \ref{thm:general1} and \ref{thm:general2}.
\begin{lem} \label{lem:1}
Let $N$ be a positive integer and $b_1,\dots,b_N \in \QQ$ such that $b_N \neq 0$. It holds that 
\begin{align*}
    [ x_1^{k_1}\cdots x_N^{k_N} ] & (1+x_1)^{b_1} (1+x_1+x_2)^{b_2} \cdots (1+x_1+\cdots + x_N)^{b_N} \\
    & = \binom{b_N}{k_N}\binom{b_{N-1}+b_N-k_{N}}{k_{N-1}} \cdots
    \binom{b_{1}+\cdots+b_N - k_N \cdots-k_{2}}{k_1}.
\end{align*}
\end{lem}  
This result contains the core identity of the present paper, since it enables the connection between the algebraic functions  $R(\textbf{x})$ of the form~\eqref{eq:genR} and hypergeometric sequences.
It can be proven in two ways: a direct approach works by multiplying the left-hand side out using the multinomial theorem, picking the needed coefficient and reducing the sum using the Chu-Vandermonde identity several times. This procedure is rather tedious and not instructive, therefore we present a combinatorially inspired proof.

\begin{proof}
Because $(1+x_1)^{b_1} \cdots (1+x_1+\cdots+x_{N-1})^{b_{N-1}}$ does not depend on $x_N$, we obtain that the left-hand side of the equation is equal to
\begin{align*}
    & [x_1^{k_1}\cdots x_N^{k_N}] \prod_{i=1}^{N-1}\bigg(1+\sum_{j=1}^i x_j\bigg)^{b_i} \cdot (1+x_1+\cdots+x_{N-1})^{b_{N}} \left(1+\frac{x_N}{1+x_1+\cdots + x_{N-1}}\right)^{b_N}\\
    &  = [x_1^{k_1}\cdots x_N^{k_N}] \prod_{i=1}^{N-1}\bigg(1+\sum_{j=1}^i x_j\bigg)^{b_i} \cdot (1+x_1+\cdots+x_{N-1})^{b_{N}} \sum_{k \geq 0} \binom{b_N}{k}\left(\frac{x_N}{1+x_1+\cdots+x_{N-1}}\right)^k \\
    & = \binom{b_N}{k_N}  [x_1^{k_1}\cdots x_{N-1}^{k_{N-1}}] (1+x_1)^{b_1} \cdots (1+x_1+\cdots+x_{N-2})^{b_{N-2}} (1+x_1+\cdots + x_{N-1})^{b_{N-1}+b_{N}-k_{N}}.
\end{align*}
Now the claim follows by iteration.
\end{proof}
We remark that the same strategy as above can be used to prove an even more general statement. Let
\[
R(\x) = \prod_{i=1}^N \bigg( 1+\sum_{j \in I_i} x_j \bigg)^{b_i},
\]
for rational numbers $b_1,\dots,b_N$, such that all variables $x_1,\dots,x_N$ appear in $R(\x)$, and sets $I_1,\dots,I_N \subseteq \{ 1,\dots, N\}$ with the property that $I_1 \cup \cdots \cup I_{n-1}  \subsetneq I_1 \cup \cdots \cup I_{n-1} \cup I_n$ for all $n=1,\dots,N$. Then, similarly to the statement in Lemma~\ref{lem:1}, the coefficient of $x^{k_1}\cdots x^{k_N}$ in $R(\x)$ is a product of binomial coefficients and the diagonal of $R(\x)$ is a generalized hypergeometric function. The notation, however, becomes quite cumbersome in this setting and no new ideas are needed; therefore we stick to the more insightful but less general case $I_n = \{ 1,\dots,n \}$.

Note that Lemma~\ref{lem:1} shares some similarities with Straub's Theorem~3.1 in~\cite{Straub14}, which provides explicit expressions
of rational power series of the form
\[ 
\Big(
(1 + x_1 + \cdots + x_{\lambda_1})
(1 +  x_{\lambda_1+1} + \cdots + x_{\lambda_1+\lambda_2})
\cdots
(1 +  x_{\lambda_1+\cdots + \lambda_{\ell-1}} + \cdots + x_N)
- 
\alpha \cdot x_1 x_2 \cdots x_N
\Big)^{-1}.
\]
In Lemma~\ref{lem:1}, we allow products of linear forms with arbitrary exponents however there is no term $\alpha\cdot x_1 x_2 \cdots x_N$, while in \cite[Theorem~3.1]{Straub14} the linear forms have disjoint variables and
appear with exponent 1. Setting $\alpha=0$ in Straub's formula also yields a product of binomial coefficients.

It is legitimate to wonder whether there is a common generalization of Lemma~\ref{lem:1} and Thm. 3.1 in~\cite{Straub14}. For instance, one may ask 
for which values of $\alpha$ is the diagonal 
\[
{\rm Diag} \left(\left( \sqrt {1- x} \left( 1-y \right) - \alpha xy \right) ^{-1} \right)
=
1+ \left( \alpha+ 1/2 \right) t+ \left( {\alpha}^{2}+2\,\alpha+3/8 \right) {t}^{2}+\cdots
 \]
 hypergeometric? For a general $\alpha$, the minimal recurrence satisfied by the coefficients of the diagonal is of order 4, for $\alpha=\pm  i/2$ it is of order 3, and it seems that the only rational value of $\alpha$ for which there exists a shorter recurrence is $\alpha=0$, in which case the diagonal is hypergeometric.  
  

\medskip
Now we want to verify a similar statement for the situation as in Theorem \ref{thm:general2}, so the case where we deal with the coefficient sequence of $(1+x_1+\cdots+x_{N-2}+2\, x_{N-1})^{b} \cdot R(x_1,\dots,x_N).$
We lay the grounds for a lemma similar to Lemma~\ref{lem:1},
by starting with a rather surprising identity.

\begin{lem}\label{lem:2}
Let $k \in \NN$ and $b \in \QQ$ arbitrary. It holds that 
\begin{align*}
    [x^k] \frac{(1+2x)^{b}}{(1+x)^{k+1}} = 4^k \binom{(b-1)/2}{k}.
\end{align*}
\end{lem}
\begin{proof}
First notice that for arbitrary $a,b$ we can compute 
\begin{align*}
[x^k] (1+2x)^{b}(1+x)^{a} = [x^k] \left( \sum_{i\geq 0} 2^i \binom{b}{i} x^i \right)\left( \sum_{j\geq 0} \binom{a}{j} x^j \right) = \sum_{j = 0}^k 2^j \binom{b}{j} \binom{a}{k-j}.
\end{align*}
So we set $a= -(k+1)$ and obtain 
\begin{align*}
    [x^k] \frac{(1+2x)^{b}}{(1+x)^{k+1}} = \sum_{j = 0}^k 2^j \binom{b}{j} \binom{-k-1}{k-j} = 2^{k}\sum_{j = 0}^k (-1)^{j} 2^{-j} \binom{b}{k-j} \binom{k+j}{k}.
\end{align*}
It remains to prove the following identity\footnote{Note that identity~\eqref{id:Gould} could alternatively be proven by using Zeilberger's creative telescoping algorithm~\cite{Zeilberger91}, or derived from identity (3.42) in \cite[p. 27]{Gould72} by setting $2n-x = b$, multiplying with $2^k$ and reverting the summation.} 
 \begin{equation}\label{id:Gould}
     \sum_{j = 0}^k (-2)^{-j} \binom{b}{k-j} \binom{k+j}{k}  = 2^k {\binom{(b-1)/2}{k}}.
\end{equation} 
To do this, we note that
\begin{align*} 
    \sum_{j = 0}^k (-1)^{j} \binom{b}{k-j} \binom{k+j}{k} u^j = \binom{b}{k} {}_2F_1 ([-k, k+1]\,;[b+1-k]\,;u),
\end{align*}
and 
\[
{}_2F_1 ([-k, k+1]\,;[b+1-k]\,;1/2) =  \frac{\Gamma((b+1-k)/2)\Gamma((b+2-k)/2)}
{\Gamma((b+1-2k)/2)\Gamma((b+2)/2)} = 2^k\frac{\binom{(b-1)/2}{k}}{\binom{b}{k}},
\]
by Kummer's identity~\cite[Eq.~3, p.~134]{Kummer1836}.
\end{proof}

The proof above  explains the special role of the coefficient $a=2$ mentioned in the introduction: it is one of the few values, along with $1$ and $-1$, for which there exists a closed form expression for the evaluation of a ${}_2F_1([\alpha,1-\alpha]\,;[\gamma]\,;u)$ at $u=1/a$.

Now we can proceed and prove the essential lemma for Theorem \ref{thm:general2}. Note that contrary to Lemma \ref{lem:1} the following statement is purely about diagonal coefficients and not for general exponents anymore. Except for the missing factor $\binom{b_{N-1}+b_N-k}{k}$ and the new two factors $4^k$ and $\binom{(b-1)/2}{k}$ the formulas are completely analogous.

\begin{lem} \label{lem:3}
Let $N$ be a positive integer and $b_1,\dots,b_{N} \in \QQ$ such that $b_{N-1}+b_N=-1$. For any $b \in \QQ$ the coefficient of $x_1^k\cdots x_N^k$ in
\begin{align*}
    (1+x_1+\cdots + x_{N-2} + 2\, x_{N-1})^{b} \cdot  (1+x_1)^{b_1} \cdots (1+x_1+\cdots + x_N)^{b_N} 
\end{align*}
is given by 
\[
    4^k\binom{(b-1)/2}{k} \binom{b_N}{k} \cdot \binom{b_{N-2} + b_{N-1} + b_N + b - 2k}{k} \cdots
    \binom{b_{1}+\cdots+b_N + b -(N-1)k}{k}.
\]
\end{lem}
\begin{proof}
By the same argument as in the proof of Lemma \ref{lem:1}, the left-hand side is equal to $\binom{b_N}{k}$ multiplied with the coefficient of $x_1^k \cdots x_{N-1}^k$ in
\begin{align*}
     (1+x_1+\cdots + x_{N-2} + 2\,x_{N-1})^{b} \cdot  \prod_{j=1}^{N-2} \left(1+\sum_{i=1}^j x_i \right)^{b_j} \cdot (1+x_1+\cdots + x_{N-2} + x_{N-1})^{b_{N-1} + b_N-k}.
\end{align*}
Because the product in the middle does not depend on $x_{N-1}$ and since we assumed $b_{N-1} + b_N = -1$, we can first compute   
\begin{align*}
    [x_{N-1}^k]  \left(1+2 \, \frac{x_{N-1}}{1+x_1+\cdots + x_{N-2}}\right)^{b} & \left(1+\frac{x_{N-1}}{1+x_1+\cdots + x_{N-2}}\right)^{-1-k} \\
    & =4^k\binom{(b-1)/2}{k}(1+x_1+\cdots+x_{N-2})^{-k},
\end{align*}
by Lemma \ref{lem:2}. Therefore we are left with 
\begin{align*}
    4^k\binom{(b-1)/2}{k}\binom{b_N}{k} \cdot [x_1^k\cdots x_{N-2}^k] \prod_{j=1}^{N-3} \left(1+\sum_{i=1}^j x_i \right)^{b_j} \cdot (1+x_1+\cdots+x_{N-2})^{b_{N-2} + b-2k-1},
\end{align*}
which is easily computed using Lemma \ref{lem:1}.
\end{proof}

Note that the requirement $b_{N-1} + b_N = -1$ comes from the $+1$ in the denominator of the left-hand side in Lemma \ref{lem:2}. Since this identity is itself surprising and does not allow for obvious generalizations, the condition on the relationship of $b_{N-1}$ and $b_N$ is necessary. 

\subsection{Proof of Theorem \ref{thm:general1}}\label{proof:thm3}
For the proof of Theorem \ref{thm:general1} we will only use Lemma \ref{lem:1} and algebraic manipulations similar to the proof of Bober's Lemma 4.1 in \cite{Bober09}.

By Lemma \ref{lem:1} we obtain the coefficient of $t^n$ for any $n \in \NN$ on the left-hand side:
\begin{align*}
    [t^n] \diag((1+x_1)^{b_1}\cdots(1+x_1+\cdots + x_N)^{b_N}) = \binom{b_N}{n} \cdots
    \binom{b_{1}+\cdots+b_N - (N-1)n}{n}.
\end{align*}
For the right-hand side we use the fact that for all $a,b$ and non-negative integers $n$ it holds
\[
    (a/b)_n((a+1)/b)_n \cdots ((a+b-1)/b)_n \cdot b^{bn} = (a)_{bn}.
\]
Then 
\[
U_k \coloneqq \prod_{i=1}^{N-k+1} (u_i^k)_n = \frac{(-b_k-\cdots -b_N)_{(N-k+1)n}}{(N-k+1)^{(N-k+1)n}},
\]
for all $k=1,\dots,N$. Similarly, 
\[
V_k \coloneqq \prod_{i=1}^{N-k} (v_i^k)_n = \frac{(-b_k-\cdots -b_N)_{(N-k)n}}{(N-k)^{(N-k)n}},
\]
for all $k=1,\dots, N-1$. Clearly $V_N \coloneqq \prod_{i=1}^{N-1} (v^N_i)_n = (n!)^{N-1}$. We deduce that 
\begin{align*}
    [t^n] {}_{M}F_{M-1}([u]\,;[v]\,;(-N)^N t) = \frac{(-N)^{nN}}{n!}\prod_{i=1}^N \frac{U_i}{V_i} = &  \frac{(-1)^{nN}}{(n!)^N}\prod_{i=1}^{N} \frac{(-b_i-\cdots -b_N)_{(N-i+1)n}}{(-b_i-\cdots -b_N)_{(N-i)n}} \\
    = & \frac{(-1)^{nN}}{(n!)^N} \prod_{i=1}^{N} (-b_i-\cdots -b_N + (N-i)n)_{n} .
\end{align*}
The claim of Theorem~\ref{thm:general1} follows from the fact that
\begin{align*}
    (-1)^n\frac{(-b_k-\cdots -b_N + (N-k)n)_{n}}{n!} & = (-1)^n \binom{-b_k-\cdots -b_N + (N-k+1)n - 1}{n} \\& = \binom{b_k + \cdots + b_N - (N-k)n}{n}.  \qquad\qquad\qquad\qquad\qquad \square 
\end{align*}

\subsection{Proof of Theorem \ref{thm:general2}}\label{proof:thm4}
The proof of Theorem \ref{thm:general2} is very similar: we will use Lemma \ref{lem:3} and the same reasoning as before. The only difference lies in the fact that because the hypergeometric function has one parameter less, we need to redefine $U_{N-1}, V_{N-1}$ and $V_{N}$. Recall that the denominator of $U_{k}$ was given by $(N-k+1)^{(N-k+1)n}$ and it cancelled with the denominator of $V_{k-1}$. In the present case, $\tilde{U}_{N-1}$ will have no denominator and therefore $2^{2n}$ from $\tilde{V}_{N-2}$ survives. This fits with the $4^k$ in the statement of Lemma \ref{lem:3} and is another indicator for the importance and essence of the constant $a=2$.

Using Lemma \ref{lem:3} we obtain the coefficient of $t^n$ for any $n \in \NN$ on the left-hand side:
\begin{align*}
    [t^n] & \diag((1+x_1+\cdots + 2\,x_{N-1} + x_N)^{b} R(x_1,\dots,x_N)) \\
    & =4^n\binom{(b-1)/2}{n} \binom{b_N}{n} \cdot \binom{b_{N-2} + b_{N-1} + b_N + b - 2n}{n} \cdots
    \binom{b_{1}+\cdots+b_N + b -(N-1)n}{n}.
\end{align*}
By the same reasoning as before, we have for all $k=1,\dots,N-2,N$
\[
\tilde{U}_k \coloneqq \prod_{i=1}^{N-k+1} (\tilde{u}_i^k)_n = \frac{(-b_k-\cdots -b_N-b)_{(N-k+1)n}}{(N-k+1)^{(N-k+1)n}},
\]
and similarly
\[
\tilde{V}_k \coloneqq \prod_{i=1}^{N-k} (\tilde{v}_i^k)_n = \frac{(-b_k-\cdots -b_N-b)_{(N-k)n}}{(N-k)^{(N-k)n}},
\]
for $k=1,\dots,N-2$. Clearly $\tilde{V}_{N-1} \coloneqq \prod_{i=1}^{N-1} (\tilde{v}^{N-1}_i)_n = (n!)^{N-1}$ and we set $\tilde{V}_N \coloneqq 1$. Moreover, this time we have $\tilde{U}_{N-1} \coloneqq (\tilde{u}^{N-1})_n = ((1-b)/2)_n$.
Altogether, we find
\begin{align*}
    [t^n] & {}_{M}F_{M-1}([\tilde{u}]\,;[\tilde{v}]\,;(-N)^N t) =  \frac{(-N)^{nN}}{n!} \prod_{i=1}^N \frac{\tilde{U}_i}{\tilde{V}_i} \\
    = &  \frac{(-1)^{nN}2^{2n} ((1-b)/2)_n}{(n!)^N}  \prod_{i=1}^{N-2} \frac{(-b_i-\cdots -b_N - b)_{(N-i+1)n}}{(-b_i-\cdots -b_N - b)_{(N-i)n}} \cdot \frac{(-b_N-b)_n}{1} \\
    = & 4^{n} \frac{(-1)^{n} ((1-b)/2)_n}{n!} \cdot \prod_{i=1}^{N-2} \frac{(-1)^n(-b_i-\cdots -b_N + (N-i)n)_{n}}{n!} \cdot \frac{(-1)^n(-b_N-b)_n}{n!}.
\end{align*}
Using the same final observation as before we conclude the proof.
\hfill $\square$

\section{Algebraicity and Hadamard grade} \label{sec:alg-grade}
      
\subsection{Algebraic cases} \label{sec:algebraic}    

We address here the following question: given $b_1,\dots,b_N \in \QQ$, $b_N\neq0$, when is the diagonal
\begin{align*}
    \diag(R(\textbf{x})) = \diag((1+x_1)^{b_1}\cdots(1+x_1+\cdots+x_N)^{b_N})
\end{align*}
an algebraic function?
\begin{coro} \label{coro}
	$\diag(R(\normalfont \textbf{x}))$ is algebraic if and only if $N=2$ and $b_2 \in \ZZ$, or $N=1$.
\end{coro}
In the proof below we will use several times the following useful
fact~\cite[Thm.~33]{Christol15}: if a generalized hypergeometric function is algebraic, then its monodromy weight is zero, that is the number of integer bottom parameters is at most equal to the number of integer top parameters.
\begin{proof}
By Theorem~\ref{thm:general1} it is sufficient to study the algebraicity of the generalized hypergeometric function $H(t)$ defined by 
   \begin{align*}
       & {}_{N(N+1)/2} F_{N(N+1)/2-1} \left([u^1,\dots,u^{N-1},-b_N]\,; [v^1,\dots,v^{N-1},1,1,\dots,1]\,; t \right),
    \end{align*}          
where $u^k$ and $v^k$, $k=1,\dots,N-1$ are defined like in \S \ref{subsec:gen-Thm1}:
\[
u^k \coloneqq \left( \frac{b}{\ell+1}, \frac{b+1}{\ell+1} , \dots , \frac{b+\ell}{\ell+1}  \right) \quad \text{ and } \quad v^{k} \coloneqq \left( \frac{b}{\ell}, \frac{b+1}{\ell} , \dots , \frac{b+\ell-1}{\ell}  \right)
\]
for $b = -(b_k+b_2+\cdots+b_N)$ and $\ell = N-k$. 

By definition, $N-1$ of the bottom parameters are ones. We claim that each tuple $u^k$ contains at most one integer and if it does contain one, then $v^k$ does as well.  From the definition of $u^k$ it follows that if some $u^k_i \in \ZZ$ then $b \in \ZZ$ and $b \equiv -i+1$ mod $\ell+1$. This shows that for any $k=1,\dots,N-1$ at most one $u_i^k \in \ZZ$. Because of the definition of $v^k$ we see that if $b \in \ZZ$ and $b \equiv i$ mod $\ell$ for some $i \in \{1,\dots,\ell\}$, then $v_i^k \in \ZZ$. This proves the claim. 

Henceforth, in order to introduce new integer parameters on the top, while not creating equally many on the bottom, it is only possible to choose $-b_N$ integer. Therefore, in order to achieve monodromy weight zero -- a necessary condition for algebraicity of $H(t)$ -- we need to have $N-1 \leq 1$. From the same argument it follows that in the case $N-1 = 1$, we need to have $-b_{N} \in \ZZ$.

Obviously for $N=1$ the diagonal is algebraic, so it remains to prove that, conversely, when $-b_2 \eqqcolon S$ is an integer and $b_1 \eqqcolon R \in \QQ$ arbitrary, then the diagonal in 
\begin{align}  \label{diag:alg}
\mathrm{Diag} \left( \frac{(1-x_1)^R}{(1-x_1-x_2)^S} \right)
  =        
 {} _{3}F_{2} \left( \left[ 
\frac{S-R}{2},\frac{S-R+1}{2},S
\right]\,;
\left[1, S-R \right]\,; 4 \, t \right). 
\end{align}  
is an algebraic function. If $R$ is an integer too, this follows from by~\cite{Polya22,Furstenberg67}. 
In the general case, one can rewrite the $_3F_2$ in~\eqref{diag:alg} as the Hadamard product
\begin{align}  \label{diag:alg-bis}      
 {} _{3}F_{2} \left( \left[ 
\frac{S-R}{2},\frac{S-R+1}{2},S
\right]\,;
\left[1, S-R \right]\,; t \right)
=
 {} _{2}F_{1} \left( \left[ 
\frac{S-R}{2},\frac{S-R+1}{2}
\right]\,;
\left[S-R \right]\,; t \right)  
\star
(1-t)^{-S}.
\end{align}  
The $_2F_1$ is algebraic as it corresponds to Case I in 
Schwarz's table~\cite{Schwarz1872}. Since $S$ is an integer,  
$(1-t)^{-S}$ is a rational function. We conclude by applying Jungen's theorem~\cite[Thm.~8]{Jungen1931}:
the Hadamard product of an algebraic and a rational function
is algebraic, see also~\cite[Prop.~6.1.11]{Stanley99}.   
\end{proof}         

\subsection{Hadamard grade} \label{sec:grade}

Recall that the \emph{Hadamard grade}~\cite{AlMe11} of a power series $S(t)$ is the least positive integer $h=h(S)$ such that $S(t)$ can be written as the Hadamard product of $h$ algebraic power series, or $\infty$ if no such product exists.
Since algebraic power series are diagonals~\cite[\S3]{Furstenberg67}, and diagonals are closed under Hadamard product~\cite[Prop.~2.6]{Christol88}, any power series with finite Hadamard grade is a diagonal~\cite[Thm.~7]{AlMe11}.    
Conversely, it is not clear whether diagonals always have finite Hadamard grade\footnote{There exist diagonals of any prescribed finite grade~\cite[Cor. 1 $\&$ 2]{RiRo14}, assuming the Rohrlich–Lang conjecture~\cite[Conj.~22]{Waldschmidt06}. If Christol's conjecture also holds, there exist diagonals of infinite grade~\cite[Prop.~1]{RiRo14}.}.


A natural question in relation with Corollary~\ref{coro} is the following: given $b_1,\ldots,b_N \in \QQ$, determine the Hadamard grade of $\diag(R(\x)) = \diag((1+x_1)^{b_1}\cdots(1+x_1+\cdots+x_N)^{b_N}) $, or at least decide if it is finite or not. We have have the following result as an application of the classification in \cite{BeHe89} (in particular case 1 in Table 8.3) and our main theorem.
\begin{coro}
    The Hadamard grade of $\diag(R(\x))$ is at most $N$. 
\end{coro}
\begin{proof}
Like in the proof of Corollary~\ref{coro}, we define 
\[
u^k \coloneqq \left( \frac{b}{\ell+1}, \frac{b+1}{\ell+1} , \dots , \frac{b+\ell}{\ell+1}  \right) \quad \text{ and } \quad v^{k} \coloneqq \left( \frac{b}{\ell}, \frac{b+1}{\ell} , \dots , \frac{b+\ell-1}{\ell}  \right)
\]
for $b = -(b_k+b_2+\cdots+b_N)$, $\ell = N-k$ and $k=1,\dots,N-1$. Then by Theorem~\ref{thm:general1} it is sufficient to study the Hadamard grade of 
\begin{align*}
   H(t) =  {}_{N(N+1)/2} F_{N(N+1)/2-1} \left([u^1,\dots,u^{N-1},-b_N]\,; [v^1,\dots,v^{N-1},1,1,\dots,1]\,; t \right).
\end{align*}      
Now notice that 
\begin{align*}
    H(t) = {}_{N}F_{N-1}([u^1]\,;[v^1]\,;t) \star \cdots \star {}_{2}F_{1}([u^{N-1}]\,;[v^{N-1}]\,;t) \star  {}_{1}F_{0}([-b_N]\,;[\ ]\,;t),
\end{align*}
and that each hypergeometric function in the Hadamard product is algebraic by~\cite[Thm.~7.1]{BeHe89}.
\end{proof}
For instance, by Theorem~\ref{thm:1}, when $S=1$, $R=1/2$, the diagonal 
\[
\mathrm{Diag} \left( \frac{(1-x_1-x_2)^R}{(1-x_1-x_2-x_3)^S} \right)
\]
is a transcendental $_2F_1$, which can be written as the Hadamard product of two algebraic functions
\[
  {} _{2}F_{1} \left( \left[ 
\frac{1}{6}, \frac{5}{6}
\right]\,;
\left[\frac12 \right]\,; t \right)  
\star  
(1-t)^{-1/2},
\]  
(the $_2F_1$ being algebraic by Schwarz's classification~\cite{Schwarz1872})
hence its Hadamard grade is 2.

\begin{figure}[]
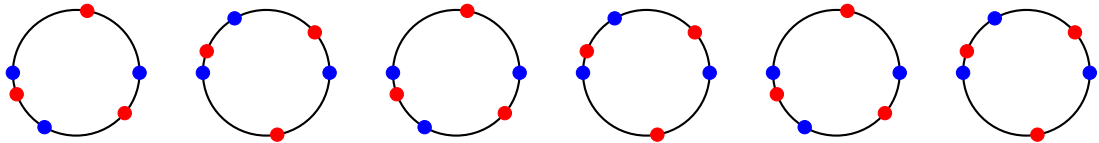
  
\begin{center} 
    \BeHe{18}{2/9}{5/9}{8/9}{1/2}{2/3}
\end{center}      
\caption{A pictorial proof of the algebraicity of the hypergeometric functions ${} _{3}F_{2} \left( \left[ 
\frac{1-R}{3},\frac{2-R}{3},\frac{3-R}{3}
\right]\,;
\left[\frac12, 1-R\right]\,; t \right) $
for $R\in\{\frac13,\frac23\}$.
There are $\varphi(18)=6$ conditions to check,
which lead to two distinct interlacing configurations.
} \label{fig:interlacing}  
\end{figure}

Similarly, the diagonals from~\eqref{eq:10} and~\eqref{eq:11} 
have Hadamard grade 2 due to the identities
 \[
  {} _{3}F_{2} \left( \left[ 
\frac{2}{9},\frac{5}{9}, \frac{8}{9}
\right]\,;
\left[1, \frac23 \right]\,; t \right)  
=
  {} _{3}F_{2} \left( \left[ 
\frac{2}{9},\frac{5}{9}, \frac{8}{9}
\right]\,;
\left[\frac12, \frac23 \right]\,; t \right)  
\star  
(1-t)^{-1/2}
\]  
and
 \[
  {} _{3}F_{2} \left( \left[ 
\frac{1}{9},\frac{4}{9}, \frac{7}{9}
\right]\,;
\left[1, \frac13 \right]\,; t \right)  
=
  {} _{3}F_{2} \left( \left[ 
\frac{1}{9},\frac{4}{9}, \frac{7}{9}
\right]\,;
\left[\frac12, \frac13 \right]\,; t \right)  
\star  
(1-t)^{-1/2}  
\]            
and to the fact that the two $_3F_2$'s on the right-hand side
are algebraic by the interlacing criterion~\cite[Thm.~4.8]{BeHe89};
see Figure~\ref{fig:interlacing} for a pictorial proof, where
red points correspond to top parameters, and blue points to bottom parameters 
(and the additional parameter 1).

More generally, the diagonal from~\eqref{eq:23-26}  
has Hadamard grade 2 due to the identity               
\begin{align*}
{} _{3}F_{2}  & \left( \left[  
\frac{1-R}{3},\frac{2-R}{3}, \frac{3-R}{3}
\right]\,;
\left[1, 1-R \right]\,; t \right)  
= \\
   & \qquad {} _{3}F_{2} \left( \left[ 
\frac{1-R}{3},\frac{2-R}{3}, \frac{3-R}{3}
\right]\,;
\left[\frac12, 1-R \right]\,; t \right)  
\star  
(1-t)^{-1/2} ,
\end{align*}              
since the $_3 F_2$ on the right-hand side is an algebraic function
for any $R\in \QQ$ (with Fig.~\ref{fig:interlacing} replaced by a similar one,
containing only interlacing blue right triangles and red equilateral triangles).   

This observation provides an alternative (and probably the shortest) proof
that the hypergeometric functions in~\eqref{eq:10}, \eqref{eq:11}
and~\eqref{eq:23-26} are diagonals of algebraic functions. However, this viewpoint does not yield such a compact diagonal representation as found in~\cite{AbKoMa20} and in our main theorem.

The same observation also quickly solves two more cases amongst the 16 cases in the list \cite[p.~58]{BBCHM12}, namely those $  {} _{3}F_{2} \left( \left[ 
N_1/9,N_2/9, N_3/9 \right]\,; \left[1, M_1/3 \right]\,; t \right) $ for which $(N_1,N_2,N_3;M_1)$ is $(1,4,7;2)$ or $(2,5,8;1)$.

Furthermore, using the interlacing criterion it is easy to see that the hypergeometric function $ {}_{3}F_{2}([1/9,4/9,7/9]\,;[a,b]\,;t)$ is algebraic if $(a,b)$ or $(b,a)$ occurs in the set  
\begin{align*}
\{     (3/4, 1/4), (2/3, 1/3), (2/3, 1/6), (1/2, 1/3), (1/2, 1/6) \}.
\end{align*}
Similarly, 
$_{3}F_{2}([2/9,5/9,8/9]\,;[a,b]\,;t)$
is algebraic if $(a,b)$ or $(b,a)$ is part of
\begin{align*}
\{    (5/6, 1/2), (5/6, 1/3), (3/4, 1/4), (2/3, 1/2), (2/3, 1/3) \}.
\end{align*}
Moreover, both 
 $   _{3}F_{2}\left(\left[\frac14, \frac38, \frac78\right]\,; \left[\frac23, \frac13 \right]\,; t\right) \ \text{and} \
    _{3}F_{2}\left(\left[\frac18, \frac34, \frac58\right]\,; \left[\frac23, \frac13 \right]\,; t\right)
$
are algebraic.      

\smallskip 
The previous analysis proves the following corollary. 
\begin{coro}   \label{coro:16cases}
	The hypergeometric function   
\begin{align*}
  {} _{3}F_{2} \left( \left[ 
A,B, C
\right]\,;
\left[1, D \right]\,; t \right)  
\end{align*}     
has Hadamard grade 2 (hence is a diagonal) for $(A,B,C;D)$ in the following set
\begin{align*}
\Big\{
(1/4, 3/8, 7/8; 1/3),  \;
1/4, 3/8, 7/8; 2/3),   \;
(1/8, 5/8, 3/4; 1/3),  \;
(1/8, 5/8, 3/4; 2/3), \\
(1/9, 4/9, 7/9; 1/2),  \;
\textcolor{red}{(1/9, 4/9, 7/9; 1/3)},  \;
(1/9, 4/9, 7/9; 1/4),  \;
(1/9, 4/9, 7/9; 1/6),    \\
\textcolor{orange}{(1/9, 4/9, 7/9; 2/3)},  \;
(1/9, 4/9, 7/9; 3/4),  \;
(2/9, 5/9, 8/9; 1/2),  \;
\textcolor{orange}{(2/9, 5/9, 8/9; 1/3)},      \\
(2/9, 5/9, 8/9; 1/4),  \;
\textcolor{red}{(2/9, 5/9, 8/9; 2/3)},  \;
(2/9, 5/9, 8/9; 3/4),  \;
(2/9, 5/9, 8/9; 5/6)
\Big \} .
\end{align*}  
\end{coro}    
 Note that the authors of~\cite{BBCHM12,BBCHM13} produced in 2011 a list of 116
potential counter-examples to Christol's conjecture; they displayed a sublist
of 18 cases in the preprint~\cite[Appendix~F]{BBCHM12}, of which they selected
3 cases that were published in~\cite[\S5.2]{BBCHM13}. As of today, to our
knowledge, the 3 cases in~\cite{BBCHM13} are still unsolved\footnote{Rivoal
and Roques proved in~\cite[Proposition 1]{RiRo14} that one of the 3 cases
in~\cite[\S5.2]{BBCHM13}, namely $_{3}F_{2}\left(\left[\frac17, \frac27,
\frac47\right]\,; \left[1, \frac12 \right]\,;t\right)$, has infinite grade assuming
the Rohrlich–Lang conjecture~\cite[Conj.~22]{Waldschmidt06}; the status of the
analogous statement for Christol's $_{3}F_{2}\left(\left[\frac19, \frac49,
\frac59\right]\,; \left[1, \frac13 \right]\,;t\right)$ is still unclear.}, while 2
of the 18 cases in~\cite{BBCHM12} have been solved in~\cite{AbKoMa20} (in red,
above) and 2 others in the current paper (in orange, above). From the list of
116 cases, only 2 were previously solved, in~\cite{AbKoMa20}.
Corollary~\ref{coro:16cases} solves 14 cases more, raising the number of solved cases to 16 (out of~116). Finally we note that another 24 cases could be resolved with the ansatz $ {} _{3}F_{2} \left( \left[ A,B, C \right]\,;\left[1, D \right]\,; t \right) = {} _{2}F_{1} \left( \left[ A,B \right]\,;\left[ r \right]\,; t \right) \star {} _{2}F_{1} \left( \left[ C,r \right]\,;\left[ D \right]\,; t \right)$ and finding $r\in\QQ$ such that both hypergeometric functions on the right-hand side are algebraic. The remaining 76 cases, however, seem to be much more difficult.

\medskip\noindent {\bf Acknowledgements.} We are grateful to Jakob Steininger,
who provided a combinatorial idea which led us to the proof of
Lemma~\ref{lem:1}, and to Herwig Hauser for his interest and engaging
questions which led to Section~\ref{sec:alg-grade}. We also thank the referees and Wadim Zudilin for their helpful comments.

\bibliographystyle{alphaabbr}
\bibliography{hyperdiag}

\end{document}